\documentclass[11pt,a4paper]{article}
\usepackage[top=3cm, bottom=3cm, left=3cm, right=3cm]{geometry}
\usepackage[latin1]{inputenc}
\usepackage{amsmath}
\usepackage{amsthm}
\usepackage{amsfonts}
\usepackage{amssymb}
\usepackage{graphicx}

\usepackage{amsmath,amsthm,amssymb,graphicx, multicol, array}
\usepackage{enumerate}
\usepackage{enumitem}
\usepackage{hyperref}

\DeclareMathOperator{\ord}{ord}

\DeclareMathOperator{\lcm}{lcm}
\newtheorem{thm}{Theorem}[section]
\newtheorem{lem}{Lemma}[section]

\newtheorem{exa}{Example}[section]

\newtheorem{dfn}{Definition}[section]

\newcommand{\N}{\mathbb{N}}
\newcommand{\Z}{\mathbb{Z}}

\usepackage{fancyhdr}
\title{The Generalized Fibonacci Primitive Roots}
\date{}
\author{N. A. Carella}
\begin{document}

\maketitle

\begin{abstract}
This note generalizes the Fibonacci primitive roots to the set of integers. An asymptotic formula for counting the number of integers with such primitive root is introduced here.\let\thefootnote\relax\footnote{\date{} 
}\\
\end{abstract}

\tableofcontents
\vskip .15 in

\section{Introduction} \label{s1}
Let \(p\geq 5\) be a prime. A primitive root \(r\neq \pm 1,s^2\) modulo \(p\) is a Fibonacci primitive root if \(r^2\equiv r+1 \bmod p \). Several basic properties of Fibonacci primitive roots were introduced in \cite{SD72}. A few other plausible generalizations are available in the literature. Fibonacci primitive roots have applications in finite field calculations, and other area of discrete mathematics, \cite{GS84}, \cite{TW14}. \\

The subset of primes with Fibonacci primitive roots is defined by 
\begin{equation}
\mathcal{P}_F=\left \{ p \in \mathbb{P}:\ord_p(r)=\varphi(p) \text{ and } r^2 \equiv r+1 \bmod p \right \} \subset \mathbb{P} . 
\end{equation}
The asymptotic formula for the counting function has the form
\begin{equation}
P_F(x)= \# \left \{ p \leq x:\ord_p(r)=\varphi(p) \text{ and } r^2 \equiv r+1 \bmod p \right \}
\end{equation}
conditional on the generalized Riemann hypothesis, was established in \cite{LH77}, and \cite{SJ88}. The complete statement of the GRH appears in Theorem \ref{thm7.1}.This note introduces a generalization to the subset of integers 
\begin{equation}\mathcal{N}_F= \left \{ n\in \mathbb{N}:\ord_n(r)=\lambda (n) \text{ and } r^2 \equiv r+1 \bmod n \right \} \subset \mathbb{N},
 \end{equation}
and provides the corresponding asymptotic counting function 
\begin{equation}
N_F(x)=\#\left\{ n\leq x:\ord_n(r)=\lambda
(n) \text{ and } r^2\equiv r+1 \bmod n \right\} .
\end{equation} \\

\begin{thm} \label{thm1.1} Assume the generalized Riemann hypothesis. Then, are infinitely many composite integers \(n\geq 1\) which have Fibonacci primitive roots. Moreover,
the number of such integers \(n\leq x\) has the asymptotic formula 
\begin{equation}
N_F(x)=\left(\frac{e^{\gamma_F-\gamma \alpha_F}}{\Gamma \left(\alpha _F\right)}+o(1)\right)\frac{x}{(\log  x)^{1-\alpha _F}}\prod _{p\in \mathcal{A}
} \left(1-\frac{1}{p^2}\right), 
\end{equation}
where  
\begin{equation}
\mathcal{A}=\left\{ p\in \mathbb{P}:\ord_p(r) =p-1,\ord_{p^2}(r)\ne p(p-1)\text{ or }
r^2 \not \equiv r+1 \bmod p^2 \right\} \nonumber, 
\end{equation}
and \(\alpha _F=0.265705\ldots \) is a constant.
\end{thm}

The constant \(\alpha _F=(27/38)\prod _{p\geq 2} (1-1/p(p-1))\) is the average density of the densities of the primes in the residue classes modulo 20 that have Fibonacci primitive roots, confer \cite{SD72} for the calculations.\\

Sections \ref{s2} to \ref{s6} provide some essential background and the required basic results. The proof of Theorem \ref{thm1.1} is settled in Section \ref{s7}. \\

\section{Some Arithmetic Functions} \label{s2}
The Euler totient function counts the number of relatively prime integers \(\varphi (n)=\#\{ k:\gcd (k,n)=1 \}\). This counting function is compactly expressed by
the analytic formula \(\varphi (n)=n\prod_{p \mid n}(1-1/p),n\in \mathbb{N} .\)\\

\begin{lem} {\normalfont (Fermat-Euler)} \label{lem2.1}If \(a\in \mathbb{Z}\) is an integer such that \(\gcd (a,n)=1,\) then \(a^{\varphi (n)}\equiv
	1 \bmod n\).
\end{lem}

The Carmichael function is basically a refinement of the Euler totient function to the finite ring \(\mathbb{Z}/n \mathbb{Z}\). 

\begin{dfn} Given an integer
\(n=p_1^{v_1}p_2^{v_2}\cdots  p_t^{v_t}\), the Carmichael function is defined by
\begin{equation}
	\lambda (n)=\lcm \left (\lambda \left(p_1^{v_1}\right),\lambda \left (p_2^{v_2}\right ) \cdots  \lambda \left (p_t^{v_t}\right ) \right )
	=\prod _{p^v  \mid \mid  \lambda (n)} p^v,
\end{equation}
where the symbol \(p^v \mid \mid n,\nu \geq 0\), denotes the maximal prime power divisor of \(n\geq 1\), and 
\begin{equation}
	\lambda
	\left(p^v\right)= \left \{
	\begin{array}{ll}
		\begin{array}{ll}
			\varphi \left(p^v\right) 
			& \text{ if } p\geq 3\text{ or } v\leq 2,  \\
			2^{v-2} 
			& \text{ if } p=2 \text{ and }v \geq 3. \\
		\end{array}
		
	\end{array} \right . 
\end{equation}
\end{dfn} 
The two functions coincide, that is, \(\varphi(n)=\lambda (n)\) if \(n=2,4,p^m,\text{ or } 2p^m,m\geq 1\). And \(\varphi \left(2^m\right)=2\lambda
\left(2^m\right)\). In a few other cases, there are some simple relationships between \(\varphi (n) \text{ and } \lambda (n)\). In fact, it seamlessly
improves the Fermat-Euler Theorem: The improvement provides the least exponent \(\lambda (n) \mid  \varphi (n)\) such that \(a^{\lambda (n)}\equiv
1 \bmod n\). \\

\begin{lem} \label{lem2.2} { \normalfont (\cite{CR10})}  Let \(n\in \mathbb{N}\) be any given integer. Then
\begin{enumerate} [font=\normalfont, label=(\roman*)]
\item The congruence \(a^{\lambda (n)}\equiv
1 \bmod n\) is satisfied by every integer \(a\geq 1\) relatively prime to \(n\), that is \(\gcd (a,n)=1\).

\item In every congruence \(x^{\lambda (n)}\equiv 1 \bmod n\), a solution \(x=u\) exists which is a primitive root \(\bmod  n\), and for any such
solution \(u\), there are \(\varphi (\lambda (n))\) primitive roots congruent to powers of \(u\).
\end{enumerate} 
\end{lem}

\begin{proof} (i) The number \(\lambda (n)\) is a multiple of every \(\lambda \left(p^v\right)=\varphi\left(p^v\right) \) such that \(p^v \mid  n\). Ergo, for any relatively prime integer \(a\geq 2\), the system of congruences 
\begin{equation}
a^{\lambda (n)}\equiv 1\bmod p_1^{v_1}, \quad a^{\lambda (n)}\equiv 1\bmod p_2^{v_2}, \quad \ldots, \quad a^{\lambda (n)}\equiv 1\bmod p_t^{v_t},
\end{equation}
where \(t=\omega (n)\) is the number of prime divisors in \(n\), is valid. 
\end{proof}

\begin{dfn}
	An integer \(u\in \mathbb{Z}\) is called a \textit{primitive root} \(\text{mod } n \) if the least exponent \(\min  \left\{ m\in \mathbb{N}:u^m\equiv 1 \bmod n \right\}=\lambda
	(n)\).
\end{dfn}

\begin{lem} \label{lem2.4}  {\normalfont (Primitive root test)} An integer $u \in \Z$ is a primitive root modulo an integer $n \in \N$ if and only if 
\begin{equation}\label{eq2-52}
u^{\lambda(n)/p} -1\not \equiv 0 \mod  n
\end{equation}
for all prime divisors $p \mid \lambda(n)$.
\end{lem}
The primitive root test is a special case of the Lucas primality test, introduced in \cite[p.\ 302]{ LE78}. A more recent version appears in \cite[Theorem 4.1.1]{CP05}, and similar sources. 

\begin{lem} \label{lem2.3}  Let $n$, and $u \in \mathbb{N}$ be integers, $\gcd(u,n)=1$. If $u$ is a primitive root modulo $p^k$ for each prime power divisor $p^k \mid n$, then, the integer $u \ne \pm 1, v^2$ is a primitive root modulo $n$.
\end{lem}

\begin{proof} Without loss in generality, let $n=pq$ with $p\geq 2$ and $q\geq 2$ primes. Let $u$ be a primitive root modulo $p$ and modulo $q$ respectively. Then
\begin{equation} \label{eq2-01}
u^{(p-1)/r}  -1\not \equiv 0 \bmod  p  \qquad  \text{  and  } \qquad u^{(q-1)/s} -1 \not \equiv 0 \bmod  q,
\end{equation}
for every prime $r \mid p-1$, and every prime $s \mid q-1$ respectively, see Lemma \ref{lem2.4}. Now, suppose that $u$ is not a primitive root modulo $n$. In particular,
\begin{equation}\label{eq2-02}
u^{\lambda(n)/t} -1\equiv 0 \bmod  n  
\end{equation}
for some prime divisor $t \mid \lambda(n)$.\\

Let $v_t(\lambda(n))$, $v_t(p-1)$, and $v_t(q-1)$ be the $t$-adic valuations of these integers. Since $\lambda(n)=\lcm(\varphi(p-1), \varphi(q-1))$, it follows that at least one of the relations
\begin{equation}\label{eq2-23}
v_t(\lambda(n))=v_t(p-1)  \qquad  \text{ or } \qquad v_t(\lambda(n))=v_t(q-1) 
\end{equation}
is valid. As consequence, at least one of the congruence equations
\begin{equation}\label{eq2-06}
u^{\lambda(n)/t} -1\equiv 0 \bmod  n  \qquad  \Longleftrightarrow \qquad u^{\lambda(n)/t}-1 \equiv \mod p
\end{equation}
or
\begin{equation}\label{eq2-08}
u^{\lambda(n)/t} -1\equiv 0 \bmod  n  \qquad  \Longleftrightarrow \qquad u^{\lambda(n)/t}-1 \equiv \mod q
\end{equation}
fails. But, this in turns, contradicts the relations in (\ref{eq2-01}) that $u$ is a primitive root modulo both $p$ and $q$. Therefore, $u$ is a primitive root modulo $n$. 
\end{proof}

\begin{exa} {\normalfont The integer $2$ is a primitive root modulo both $p=37$ and $q=61$. Let $n=p\cdot q=37 \cdot 61$, $\varphi(p-1)=2^2\cdot 3^2$, $\varphi(q-1)=2^2\cdot 3 \cdot 5$, and $\lambda(n)=\lcm(\varphi(p-1), \varphi(q-1))=2^2\cdot 3^2 \cdot 5$. The corresponding congruences and $t$-adic valuations of these integers are these.
\begin{itemize}
\item For $t=2$, the valuations are: $v_2(\lambda(n))=v_2(p-1)=v_2(q-1)=2$. The assumption that 2 is not a primitive root modulo $n$ is not valid: 
\begin{equation} \label{eq2-72}
u^{\lambda(n)/2} -1\equiv 0 \bmod  n  
\end{equation}
fails because at least one 
\begin{equation} \label{eq2-76}
2^{\lambda(n)/2}-1 \not \equiv 0\mod p  \qquad  \text{ or } \qquad 2^{\lambda(n)/2}-1 \not \equiv 0\mod q
\end{equation}
contradicts it.
\item For $t=3$, the valuations are: $v_3(\lambda(n))=v_3(p-1)=2$, and $v_3(q-1)=1$. The assumption that 2 is not a primitive root modulo $n$ is not valid: 
\begin{equation}\label{eq2-78}
u^{\lambda(n)/3} -1\equiv 0 \bmod  n  
\end{equation}
fails because at least one 
\begin{equation}\label{eq2-79}
2^{\lambda(n)/3}-1 \not \equiv 0\mod p  \qquad  \text{ or } \qquad 2^{\lambda(n)/3}-1 \equiv 0 \mod q
\end{equation}
contradicts it.
\item For $t=5$, the valuations are: $v_5(\lambda(n))=v_5(q-1)=1$, and $v_5(p-1)=0$. The assumption that 2 is not a primitive root modulo $n$ is not valid: 
\begin{equation}\label{eq2-78}
u^{\lambda(n)/5} -1\equiv 0 \bmod  n  
\end{equation}
fails because at least one 
\begin{equation}\label{eq2-79}
2^{\lambda(n)/5}-1  \equiv 0\mod p  \qquad  \text{ or } \qquad 2^{\lambda(n)/5}-1 \not \equiv 0 \mod q
\end{equation}
contradicts it.
\end{itemize}
Since the congruence (\ref{eq2-72}) fails for every prime divisor $t=2,3,5$ of $\lambda(n)=2^2\cdot 3^2 \cdot 5$, it implies that $2$ is primitive root modulo $n=37 \cdot 61$.
}
\end{exa}

\section{ Characteristic Function For Fibonacci Primitive Roots} \label{s3}
The symbol \(\text{ord}_{p^k}(u)\) denotes the order of an element \(r\in \left(\mathbb{Z}\left/p^k\right.\mathbb{Z}\right)^{\times }\) in the multiplicative
group of the integers modulo \(p^k\). The order satisfies the divisibility condition \(\text{ord}_{p^k}(r) \mid \lambda (n)\), and
primitive roots have maximal orders \(\text{ord}_{p^k}(r)=\lambda (n)\). The basic properties of primitive root are explicated in \cite{AT76}, \cite{RH94}, et cetera. The characteristic function \(f:\mathbb{N}\longrightarrow \{ 0, 1 \}\) of a fixed primitive root \(r\) in the finite ring \(\mathbb{Z}\left/p^k\right.\mathbb{Z}\),
the integers modulo \(p^k\), is determined here. 
\\

\begin{lem} \label{lem3.1}Let \(p^k,k\geq 1\), be a prime power, and let \(r\in \mathbb{Z}\) be an integer such that \(\gcd \left(r,p^k\right)=1\).
Then
\begin{enumerate} [font=\normalfont, label=(\roman*)]
\item  The characteristic \(f\) function of the primitive root \((r \bmod p^k)\) is given by \\ 
\begin{equation}
f\left(p^k\right) =\left \{
\begin{array}{ll}
 \begin{array}{ll}
 0 & \text{    }\text{ if } p \leq 3, k\leq 1, \\
 1 & \text{    }\text{ if } \ord_{p^k}(r)= p^{k-1}(p-1) \text{ and }
r^2 \equiv r+1 \mod p^k , \\
 &\text{ for any }p \geq 5, k \geq 1, \\
 0 & \text{    }\text{ if } \ord_{p^k}(r)\neq p^{k-1}(p-1), \text{ and } p \geq 5, k\geq 2. \\
\end{array}
\end{array} \right .
\end{equation}

\item  The function \(f\) is multiplicative, but not completely multiplicative since 
\item $f(p q) =f(p)f(q), \gcd (p,q)=1$, 
\item $f\left(p^2\right) \neq f(p)f(p), \text{     if }\ord_{p^2}(r)\neq p(p-1)$.
\end{enumerate}
\end{lem}

\begin{proof}  The congruences \(x^2-x-1\equiv 0 \text{ mod }2\) and \(x^2-x-1\equiv 0 \bmod 3\) have no roots, so \(f\left(2^k\right)=0\)
and \(f\left(3^k\right)=0\) for all \(k\geq 1.\) Moreover, for primes \(p\geq 5 \text{ and } k\geq 1,\) function has the value \(f\left(p^k\right)=1\)
if and only if the element \(r\in \left(\mathbb{Z}\left/p^k\right.\mathbb{Z}\right)^{\times }\) is a Fibonacci primitive root modulo \(p^k\). Otherwise,
it vanishes: \(f\left(p^k\right)=0\). The completely multiplicative property fails for \(p=5\). Specifically, \(0=f\left(5^2\right)\neq f(5)f(5)=1\).
\end{proof}

Observe that the conditions \(\text{ord}_p(r)=p-1 \text{ and } \text{ord}_{p^2}(r)\neq p(p-1)\) imply that the integer \(r\neq \pm 1,s^2\) cannot be
extended to a primitive root \(\bmod p^k,k\geq 2\). But that the condition \(\ord_{p^2}(r)=p(p-1)\) implies that the integer \(r\) can
be extended to a primitive root \(\bmod  p^k,k\geq 2\).

\section{Wirsing Formula} \label{s4}
This formula provides decompositions of some summatory multiplicative functions as products over the primes supports of the functions. This technique works well with certain multiplicative functions, which have supports on subsets of primes numbers of nonzero densities.

\begin{lem} {\normalfont (\cite[p. 71]{WE61}) }  \label{lem4.1}Suppose that \(f:\mathbb{N}\longrightarrow \mathbb{C}\) is a multiplicative function with the following properties.
\begin{enumerate} [font=\normalfont, label=(\roman*)]
\item $f(n) \geq 0$ for all integers $n\in \mathbb{N}$.
\item $f\left(p^k\right)\leq c^k$ for all integers $k\in \mathbb{N}$, and $c<2$ constant. 
\item $\displaystyle \sum _{p\leq x} f(p)=(\tau +o(1)) x/\log  x$, where $\tau >0$ is a constant as $x \longrightarrow  \infty$. 
\end{enumerate}
Then 
\begin{equation}
\sum _{n\leq x} f(n)=\left(\frac{1}{e^{\gamma \tau }\Gamma (\tau )}+o(1)\right)\frac{x}{\log  x}\prod _{p\leq x} \left(1+\frac{f(p)}{p}+\frac{f\left(p^2\right)}{p^2}+\cdots
\right) .
\end{equation}
\end{lem}

The gamma function appearing in the above formula is defined by \(\Gamma (s)=\int _0^{\infty }t^{s-1}e^{-s t}d t, s\in \mathbb{C}\). The intricate
proof of Wirsing formula appears in \cite{WE61}. It is also assembled in various papers, such as \cite{HA87}, \cite[p. 195]{PA88}, and discussed in \cite[p. 70]{MV07}, \cite[p. 308]{TG15}. Various applications are provided in \cite{MP11}, \cite{PS03}, \cite{WK75}, et alii. \\

\section{Harmonic Sums And Products Over Fibonacci Primes} \label{s5}
The subset of primes \(\mathcal{P}_F=\left\{ p\in \mathbb{P}:\ord(r)=\varphi (p) \text{ and } r^2\equiv r+1 \mod p \right\}\subset \mathbb{P}\) consists of all the primes with Fibonacci primitive roots \(r\in \mathbb{Z}\). By Theorem \ref{thm7.1}, it has nonzero density \(\alpha _F=\delta \left(\mathcal{P}_F\right)>0\). The real number \(\alpha _F>0\) is a rational multiple of the well known Artin constant, see \cite{SD72}, and \cite{LH77} for the derivation. The proof of the next result is based on standard analytic number theory methods in the literature, refer to \cite[Lemma 4]{PS03}.

\begin{lem} \label{lem5.1} Assume the generalized Riemann hypothesis, and let \(x\geq 1\) be a large number. Then, there exists a pair
of constants \(\beta _F>0,\text{     and     } \gamma _u>0\) such that 
\begin{enumerate} [font=\normalfont, label=(\roman*)]
\item $ \displaystyle \sum _{\substack{p\leq x\\ p\in \mathcal{P}_F}} \frac{1}{p}=\alpha _F\log \log
x+\beta _F+ O\left(\frac{\log \log x}{\log^2  x}\right) .$ \\
\item $ \displaystyle 
\sum _{\substack{p\leq x\\ p\in \mathcal{P}_F}} \frac{\log p}{p-1}=\alpha _F \log
x-\gamma _F+ O\left(\frac{\log \log x}{\log^2 x}\right) $.
\end{enumerate}
\end{lem}

\begin{proof} (i) Let \(\pi _F(x)=\#\left\{ p\leq x:\ord(r)=\varphi (p) \text{ and } r^2\equiv r+1 \text{ mod }p \right\}\) be the counting measure of the corresponding subset of primes \(\mathcal{P}_F\). To estimate the asymptotic order of the prime harmonic
sum, use the Stieltjes integral representation:
\begin{equation}
\sum _{\substack{p\leq x\\ p\in \mathcal{P}_F}} \frac{1}{p} =\int_{x_0}^{x}\frac{1}{t}d \pi _F(t) =\frac{\pi _F(x)}{x}+c_(x_0)
+\int_{x_0}^{x}\frac{\pi _u(t)}{t^2}d t,
\end{equation} \\
where \(x_0>0\) is a constant. By the GRH, $\pi_F(x)=\alpha
_F\pi (x)=\alpha_Fx/ \log x+O\left(\log \log x /\log^2  x\right)$, see Theorem 7.1. This yields
\begin{eqnarray}
\int_{x_0}^{x}\frac{1}{t}d \pi _F(t) 
&=&\frac{\alpha _F}{\log  x}+O\left(\frac{\log \log x}{\log ^2(x)}\right)+c_0(x_0) \nonumber\\
&\qquad&+\alpha_F\int_{x_0}^{x}\left(\frac{1}{t \log  t}+O\left(\frac{\log \log t}{t \log ^2(t)}\right)\right)d t \\
&=&\alpha _F \log \log x-\log \log x_0+c_0\left(x_0\right)+O\left(\frac{\log \log x}{\log^2  x}\right)  \nonumber,
\end{eqnarray} \\
where \(\beta _F=-\text{loglog} x_0+c_0\left(x_0\right)\) is a generalized Mertens constant. The statement (ii) follows from statement (i) and partial summation.
\end{proof}

A generalized Mertens constant \(\beta _F\) and a generalized Euler constant \(\gamma _F\) have other equivalent definitions such as \\
\begin{equation}  \label{el330}
\beta _F=\lim_{x \rightarrow
 \infty } \left(\sum _{\substack{p\leq x\\ p\in \mathcal{P}_F}} \frac{1}{p}-\alpha _u\log \log x\right)\text{              }\text{   and   }\text{            }\beta
_u=\gamma _u-\text{  }\sum _{p\in \mathcal{P}_F,} \sum _{k\geq 2} \frac{1}{k p^k} ,
\end{equation}\\
respectively. These constants satisfy \(\beta _F=\beta _1\alpha _F\text{ and } \gamma _F=\gamma \alpha _F\). If the density \(\alpha _F=1\), these definitions reduce to the usual Euler constant and the Mertens constant, which are defined by the limits 
\begin{equation}
\gamma =\lim_{x \rightarrow  \infty }
\left ( \sum _{ p\leq x } \frac{\log  p}{p-1}-\log  x\right) \text{ and } \beta _1=\lim_{x \rightarrow  \infty } \left(\sum _{ p\leq x } \frac{1}{p}-\log \log x \right),
\end{equation}
or some other equivalent definitions, respectively. Moreover, the linear independence relation in (\ref{el330}) becomes \(\text{  }\beta =\gamma -\text{ 
}\sum _{p\geq 2} \sum _{k\geq 2} \left(k p^k\right) ^{-1}\), see \cite[Theorem 427] {HW08}.
\\

A numerical experiment utilizing a small subset of primes with Fibonacci primitive roots gives the followings approximate values: 
\begin{enumerate}
\item $ \displaystyle \alpha _F=\frac{27}{38}\prod
_{p\geq 2} \left(1-\frac{1}{p(p-1)}\right)=0.265705484288843681890137 \ldots,$ \\
\item $ \displaystyle \beta _F \approx \sum _{ \substack{p\leq 1301\\ p\in \mathcal{P}_F}} \frac{1}{p} -\alpha
_F\log \log x=0.05020530308647012230491, $    and \\
\item $ \displaystyle  \gamma _F\approx \sum _{ \substack{p\leq 1301\\ p\in \mathcal{P}_F}} \frac{\log
 p}{p-1} -\alpha _F  \log  x=0.0221594862523476326826286 \ldots . $
\end{enumerate}

\begin{lem}   \label{lem5.2} Assume the generalized Riemann hypothesis, and let \(x\geq 1\) be a large number. Then, there exists a pair of constants $\gamma _F>0$ and $\nu _F >0$ such that 
\begin{enumerate} [font=\normalfont, label=(\roman*)]
\item $ \displaystyle  \prod _{\substack{p\leq x\\ p\in \mathcal{P}_F}} \left(1-\frac{1}{p}\right)^{-1}=e^{\gamma
_F} \log (x)^{\alpha _F}+O\left(\frac{\log \log x}{\log^2  x}\right) .$ \\
\item $  \displaystyle \prod _{\substack{p\leq x p\in \mathcal{P}_F}} \left(1+\frac{1}{p}\right)=e^{\gamma
_F}\prod _{p\in \mathcal{P}_F} \left(1-p^{-2}\right) \log (x)^{\alpha _F}+O\left(\frac{\log \log x}{\log^2  x}\right) . $ 
\item $ \displaystyle \prod
_{\substack{p\leq x\\ p\in \mathcal{P}_F}} \left(1-\frac{\log  p}{p-1}\right)^{-1}=e^{\nu _F-\gamma _F}x^{\alpha _F}+O\left(\frac{x^{\alpha _F}\log \log
x}{\log^2  x}\right) . $
\end{enumerate}
\end{lem}

\begin{proof} (i) Express the logarithm of the product as \\
\begin{eqnarray}
\sum _{\substack{p\leq x\\ p\in \mathcal{P}_F}} \log \left(1-\frac{1}{p}\right)^{-1}&=&\sum
_{p\leq x, p\in \mathcal{P}_F,} \sum _{k\geq 1} \frac{1}{k p^k}\\
&=&\sum _{\substack{p\leq x\\ p\in \mathcal{P}_F}} \frac{1}{p} +\sum _{p\leq x, p\in \mathcal{P}_F,}
\sum _{k\geq 2} \frac{1}{k p^k} .
\end{eqnarray}   \\

Apply Lemma 5.1 and relation (\ref{el330}) to complete the verification. For (ii) and (iii), use similar methods as in the first one.  
\end{proof}

The constant \(\nu _F>0\) is defined by the double power series (an approximate numerical value for \(p\leq x=1301\) is shown): 
\begin{equation}
\nu _2=\sum _{p\in \mathcal{P}_F, } \sum _{k\geq 2} \frac{1}{k}\left(\frac{\log  p}{ p-1}\right)^k\approx 0.188622600886988493134287\text{... }.
\end{equation}\\

\section{Density Correction Factor} \label{s6}
The subsets of primes \(\mathcal{P}_F=\left\{ p\in \mathbb{P}:\ord(r)=p-1 \text{  and  } r^2=r+1 \right\}\) has an important partition as a disjoint
union \(\mathcal{P}_F=\mathcal{A}\cup \mathcal{B}\), where 
\begin{equation} \mathcal{A}=\left \{ p\in \mathbb{P}:\text{ord}(r)=p-1,\ord_{p^k} r \ne
p^{k-1}(p-1),\text{ or } r^2\not \equiv r+ 1 \text{ mod }p^2 \right \} 
\end{equation}   and   
\begin{equation}\mathcal{B}=\left \{ p\in \mathbb{P}:\ord(r)=p-1,\ord_{p^k}
r=p^{k-1}(p-1), \text{ and } r^2\equiv r+1 \text{ mod }p^2 \right\}.
\end{equation}   

This partition has a role in the determination of the density of the subset of integers \(\mathcal{N}_F=\left\{ n\in \mathbb{N}:\text{ord}(r)=\lambda
(n) \text{ and } r^2\equiv r+1 \text{ mod }n \right\}\) with Fibonacci primitive roots \(r\in \mathbb{Z}\). The finite prime product occurring in the proof
of Theorem \ref{thm1.1} splits into two finite products over these subsets of primes. These are used to reformulate an equivalent expression suitable for this
analysis:

\begin{eqnarray} \label{el550}
P(x)&=&\prod_{\substack{p^k\leq
x,\\ \ord_p(r)=p-1, \\ \ord_{p^2}(r)\neq p(p-1) \text{ or }r^2 \not \equiv r+1 \bmod  p^2}} \left(1+\frac{1}{p}\right) \nonumber \\
&\times &\prod _{\substack{p^k\leq x,\\ \ord_{p^2}(r)=p(p-1) \text{ and }r^2 \equiv r+1 \bmod  p^2}} \left(1+\frac{1}{p}+\frac{1}{p^2}+\cdots \right)  \nonumber \\
&=&\prod_{\substack{p^k\leq
x,\\ \ord_p(r)=p-1, \\ \ord_{p^2}(r)\neq p(p-1) \text{ or }r^2 \not \equiv r+1 \bmod  p^2}} \left(1-\frac{1}{p^2}\right)  \\
&\times& \prod _{\substack{p^k\leq x,\\ \text{ord}_{p^2}(r)=p(p-1) \text{ and }r^2 \equiv r+1 \bmod  p^2}} \left(1-\frac{1}{p} \right)^{-1}   \nonumber\\
&=&\prod _{\substack{p \leq x\\ p\in \mathcal{A}}} \left(1-\frac{1}{p^2}\right)\prod _{\substack{p\leq x\\p\in \mathcal{P}_F }} \left(1-\frac{1}{p}\right)^{-1}  \nonumber \\ 
&=&\prod _{p\in \mathcal{A}} \left(1-\frac{1}{p^2}\right)\prod _{\substack{p\leq x\\ p\in \mathcal{P}_F} } \left(1-\frac{1}{p}\right)^{-1} +O\left(\frac{\log x}{x}\right) \nonumber , 
\end{eqnarray} 
where the convergent partial product is replaced with the approximation 
\begin{equation}
\prod _{\substack{p\leq x\\ p\in \mathcal{A}} }  \left(1-\frac{1}{p^2}\right)=\prod
_{p\in \mathcal{A} }  \left(1-\frac{1}{p^2}\right)+O\left(\frac{1}{x}\right) .
\end{equation} 

The product \(\prod _{p\in \mathcal{A} }  \left(1-p^{-2}\right)\) reduces the density to compensate for those primes for which the Fibonacii primitive
root \(r \bmod p\) that cannot be extended to a Fibonacci primitive root \(r \bmod p^2\). This seems to be a density correction factor similar
to case for primitive roots over the prime numbers. The correction required for certain densities of primes with respect to fixed primitive roots
over the primes was discovered by the Lehmers, see \cite{SP03}.

\section{The Proof Of The Theorem} \label{s7}
The algebraic constraint \(r^2=r+1\) restricts the Fibonacci primitive roots \(r\) modulo a prime \(p\geq 5\) to the algebraic integers \(\left.\left(1\pm \sqrt{5}\right)\right/2\) in the finite field \(\mathbb{F}_p\). This in turns forces the corresponding primes to the residue classes \(p\equiv \pm 1 \text{ mod }10\), the analysis appears in \cite{SD72}. A conditional proof for the conjecture of Shank was achieved in \cite{LH77} and \cite{SJ88}.

\begin{thm} {\normalfont (\cite{LH77}) } \label{thm7.1}If the extended Riemann hypothesis hold for the Dedekind zeta function over Galois fields of the
type \(\mathbb{Q}\left(\sqrt[d]{r},\sqrt[n]{1}\right)\), where \(n\) is a squarefree integer, and \(d|n.\) Then, the density of primes which have
Fibonacci primitive roots is \(\alpha _F=(27/38)\prod _{p\geq 2} (1-1/p(p-1))\), and the number of such primes is
\begin{equation}
P_F(x)=\alpha
_F\frac{x}{\log  x}+O\left(\frac{x \log \log x}{\log ^2 x}\right)\text{        }\text{ as  } x\longrightarrow \infty  .
\end{equation}
\end{thm}

In terms of the subsets \(\mathcal{A} \text{ and } \mathcal{B}\), see Section \ref{s6}, the characteristic function \(f(n)\) of Fibonacci primitive
roots has the simpler description 
\begin{equation}
f\left(p^k\right) =
\left \{ \begin{array}{ll}
 1 & \text{    }\text{ if } p\in \mathcal{A}, \text{ and }\text{  }k=1, \\
 0 & \text{    }\text{ if } p\in \mathcal{A}, \text{ and }\text{  }k\geq 2, \\
\end{array} \right .
\text{          }\text{ and }\text{            }f\left(p^k\right) = \left \{ 
\begin{array}{ll}
 1 & \text{    }\text{ if } p\in \mathcal{B}, \\
 0 & \text{    }\text{otherwise}. \\
\end{array} \right .
\end{equation} 

\begin{proof} By Theorem \ref{thm7.1}, the density \(\alpha _F=\delta \left(\mathcal{P}_F\right)>0\) of the subset of primes
\(\mathcal{P}_F\) is nonzero. Put \(\tau =\alpha _F\) in Wirsing formula, Lemma \ref{lem4.1}, and replace the characteristic function \(f(n)\) of Fibonacci
primitive roots in the finite ring \(\mathbb{Z}\left/p^k\right. \mathbb{Z},k\geq 1\), see Lemma \ref{lem3.1}, to produce

\begin{eqnarray}
\sum _{n\leq x} f(n)&=&\left(\frac{1}{e^{\gamma
\tau }\Gamma (\tau )}+o(1)\right)\frac{x}{\log  x}\prod _{p\leq x} \left(1+\frac{f(p)}{p}+\frac{f\left(p^2\right)}{p^2}+\cdots \right) \text{   
             } \nonumber\\
&=&\left(\frac{1}{e^{\gamma \alpha _F}\Gamma \left(\alpha _F\right)}+o(1)\right)\frac{x}{\log  x}    \nonumber\\
&\times& \prod_{\substack{p^k\leq
x,\\ \text{ord}_p(r)=p-1, \\ \text{ord}_{p^2}(r)\neq p(p-1) \text{ or }r^2 \not \equiv r+1 \text{mod } p^2}} \left(1+\frac{1}{p}\right) \\
& \times& \prod _{\substack{p^k\leq x,\\ \text{ord}_{p^2}(r)=p(p-1) \text{ and }r^2 \equiv r+1 \text{mod } p^2}} \left(1+\frac{1}{p}+\frac{1}{p^2}+\cdots \right)   \nonumber.
\end{eqnarray}

Replacing the equivalent product, see (\ref{el550}) in Section \ref{s5}, and using Lemma \ref{lem5.2}, yield\\
 \begin{eqnarray}
\sum _{n\leq x} f(n)&=&\left(\frac{1}{e^{\gamma \alpha
_F}\Gamma \left(\alpha _F\right)}+o(1)\right)\frac{x}{\log  x} \nonumber \\
&=&\prod _{p\in \mathcal{A} }  \left(1-\frac{1}{p^2}\right)\prod _{\substack{p\leq x\\ p\in \mathcal{P}_F
}} \left(1-\frac{1}{p}\right)^{-1} \\
&=&\left(\frac{e^{\gamma _F-\gamma \alpha _F}}{\Gamma \left(\alpha _F\right)}+o(1)\right)\frac{x}{(\log
 x)^{1-\alpha _F}}\prod _{p\in \mathcal{A} }  \left(1-\frac{1}{p^2}\right)   \nonumber
\end{eqnarray}
where \(\alpha _F>0\) is a constant, and \(\gamma _F\) is a generalized Euler constant, see Lemma 6 for detail. 
\end{proof}

As stated before, the constant \(\alpha _F=(27/38)\prod _{p\geq 2} (1-1/p(p-1))=0.265705\text{...}\) is the average density of the densities of the primes in the residue
classes mod 20 that have Fibonacci primitive roots, confer \cite{SD72} for the calculations.\\

The Fibonacci primitive root problem is a special case of a more general problem that investigates the existence of infinitely many primes such that the polynomial congruence $f(n) \equiv 0 \bmod p$ has primitive root solutions. The algebraic and analytic analysis for these more general polynomials are expected to be interesting.


\end{document}